\documentclass[a4paper,reqno]{amsart}

\usepackage[all]{xy}
\usepackage{color}
\usepackage{mathrsfs}
\usepackage{hyperref}
\usepackage{enumerate}
\usepackage{amsmath}
\usepackage{amsthm}
\usepackage{amssymb}
\usepackage{amscd}
\usepackage{graphicx}
\usepackage{epsfig}
\usepackage[english]{babel}
\usepackage[stable]{footmisc}

\newtheorem{proposition}{Proposition}
\newtheorem{lemma}{Lemma}
\newtheorem{theorem}{Theorem}
\newtheorem{corollary}{Corollary}
\newtheorem{definition}{Definition}

\newtheorem{remark}{Remark}

\begin{document}

\author{D. Mart\'inez Torres}

\email{dfmtorres@gmail.com}

\address{Centro de An\'{a}lise Matem\'{a}tica, Geometria e Sistemas Din\^{a}micos,
Departamento de Matem\'atica, Instituto Superior T\'ecnico, Av. Rovisco Pais, 1049-001
Lisboa, Portugal}

% \thanks{Partially supported by Fundaci\'on Universidad Carlos III,
%  the European Human Potential Program through EDGE, HPRN-CT-2000-00101,
%   the Galileo Galilei postdoctoral program at Pisa University, and 
% research projects MTM2004-07090-C03 of the Spanish Ministry of
% Science and Technology, Symmetries and deformations in geometry from the
%  NWO Council for Physical sciences, Geometry and Quantum theory from the 
% NWO Dutch Research Council, and Funda\c{c}\~{a}o para a Ci\^{e}ncia e a Tecnologia (FCT /
% Portugal).}

\title{Generic linear systems for projective CR manifolds}

\begin{abstract}  For compact CR manifolds of hypersurface type which
embed in complex projective space, we show that for all k large enough
there exist linear systems of ${\mathcal{O}}(k)$ which when
restricted to the CR manifold are generic in a suitable sense. These systems are constructed using approximately
holomorphic geometry.
\end{abstract}
% 
% \vspace{.2cm}
\maketitle
% MSC: 32V10 (primary); 58A20, 53A20, 51N15 (secondary).
% 
% Keywords: CR structure; CR jet bundle; generic linear system; Lefschetz pencil; dual geometry; approximately holomorphic geometry.
% 

\section{Introduction}\label{sec:sec1}

Approximately holomorphic geometry can be applied to compact symplectic
 manifolds endowed with a compatible almost complex structure  to obtain ``generic'' approximately holomorphic maps to complex
  projective spaces \cite{Au01}. These maps can be understood as 
   analogs in symplectic geometry of generic linear systems for Hodge manifolds. The analogy is valid 
   for constructions only using (pseudo-holomorphic)  1 and 2-jets
     \cite{Do98,Au00}; it breaks down for higher order jets  due to the difficulty
      of developing a theory of normal forms in the absence of an integrable almost
       complex structure. However if the symplectic structure comes from a Hodge one, then according to
        section 7 in \cite{Do96} it is possible to adjust the constructions of
         approximately holomorphic geometry so that the outcome are holomorphic maps, giving
          thus a new construction of generic linear systems.

Applications of approximately holomorphic theory to symplectic manifolds include
 the construction of plenty of symplectic submanifolds with control on their topology
 \cite{Do96,Au97,MPS02}, construction of symplectic invariants  \cite{AK00},
  and a proof -via Lefschetz pencil structures \cite{Do98} and making little
  use of elliptic theory- of the existence of symplectic curves realizing the
  canonical class of a 4-dimensional symplectic manifold \cite{DS00}.

A 2-calibrated manifold is defined as a triple
$(M^{2n+1},D^{2n},\omega)$,
 where $D$ is a codimension one distribution and $\omega$ a closed 2-form maximally
  non-degenerate over $D$ (i.e. making $D$ into a symplectic distribution).  The most relevant 
classes of 2-calibrated manifolds
   are related to the two extreme  behaviors of the distribution: when $D$ is
    integrable we speak of 2-calibrated foliations, and when $D$ is maximally
     non-integrable 2-calibrated structures include contact structures.

Approximately holomorphic geometry can be adapted to compact 2-calibrated
 manifolds to obtain  ``generic'' approximately holomorphic maps to complex
  projective spaces \cite{Ma04}. Very much as in the symplectic case, due the difficulty of developing
a theory of normal forms  interesting applications only use (pseudo-holomorphic) 1-jets. They include the existence of contact
    submanifolds  whose Poincar\'e dual realizes any determinantal class
     \cite{Ma04}, the construction of open book decompositions compatible
      with a contact structure \cite{GM01,Gi02}, and a description of
       the leaf spaces of 2-calibrated foliations \cite{Ma05}. 

When the auxiliary compatible almost complex structure on $D$ needed to
 develop the approximately holomorphic theory  is integrable,  we also 
  have a CR structure (of hypersurface type). In this paper we are interested
   in analyzing under which conditions the approximately holomorphic
    techniques can be refined to yield CR constructions, and which
      are the applications that can be obtained. To this extent our main result is 
the following theorem (the reader is referred to section \ref{sec:sec2} for a full explanation of its statement):

\begin{theorem}\label{thm:gener0} Let $(M^{2n+1},{\mathcal{F}},J)$ be a closed Levi-flat
 CR manifold of hypersurface type endowed with a positive CR line bundle. Fix
 h,r $\in \mathbb N$, r $\leq$ h-2. Then for any integer m there exists
 $\phi\colon M\backslash B\rightarrow {\mathbb CP}^m$ an  r-generic map.
 More precisely, we obtain the following:
\begin{enumerate}
\item A Levi-flat CR submanifold $B$ of real codimension 2m+2 and class $C^h$.
\item A CR map $\phi\colon M\backslash B \rightarrow {\mathbb CP}^m$ such that
  for each leaf $F$ of $\mathcal{F}$, the holomorphic  r-jet of  $\phi_{\mid F}$ 
   is transverse to the corresponding Thom-Boardman stratification of
    the bundle of holomorphic r-jets of holomorphic maps from $F$ to ${\mathbb CP}^m$. These bundles fit into a bundle
    of class $C^{h-r}$ -the bundle of CR $r$-jets of CR maps from $M$ to ${\mathbb CP}^m$- and the same holds
    for the strata of the Thom-Boardman stratifications. The CR r-jet
   of $\phi$  is leafwise transverse to this stratification by Levi-flat manifolds
      of class $C^{r-h}$. Therefore, the pullback of each stratum is a
       Levi-flat submanifold of the expected codimension and of class $C^{r-h}$.
\end{enumerate}
\end{theorem}

For foliated manifolds there is an obvious notion of leafwise genericity, and the existence of
 leafwise r-generic maps  (holomorphic, smooth) is in general obstructed; such a map would give 
rise to a stratification of the manifold transverse to the foliation whose existence might be
 in general not possible for topological reasons.  A simplified interpretation of theorem 
\ref{thm:gener0} is that for a certain class of manifolds foliated by complex leaves, one
 can find -away from suitable submanifolds- leafwise holomorphic r-generic maps. 

There are related results in the literature addressing the embedding problem for CR
manifolds into projective space \cite{BdM75,OS00,MY04} (thus only involving CR
1-jets). More generally E. Ghys \cite{Gh99} and B. Deroin \cite{De08} have shown
that compact manifolds laminated by complex leaves -subject to some additional
conditions- do admit enough (meromorphic) functions, so that leafwise immersions
in projective spaces are possible. To the best of our knowledge our r-genericity
result, for arbitrary r, is new.

The structure of the paper is as follows: in section \ref{sec:sec2}
 we recall the definitions  and results from CR geometry needed
  to state the results of this paper. For complex projective manifolds  classical generic linear systems
are constructed using dual geometry.  In section  \ref{sec:sec3} we outline the difficulties
  to develop  a dual geometry for arbitrary projective CR manifolds
   (i.e. those which have a CR embedding in projective space); we mention very briefly a family of 
CR manifolds for which the classical approach extends.  Finally, in section \ref{sec:sec4} we sketch
 how to adapt the constructions
 of  approximately holomorphic geometry for projective CR manifolds,
  thus proving the results stated in section \ref{sec:sec2}.

\section{Definitions and statements of the results}\label{sec:sec2}

\begin{definition} A CR manifold (always of hypersurface type for us)
 is a triple $(M^{2n+1},D^{2n},J)$, where $M$ is a manifold and $D$ is a codimension 
one distribution of $TM$ endowed with an almost complex structure $J$, such that either of the eigen-bundles
  $D^{*1,0}, D^{*0,1}$ of $D^*\otimes_{\mathbb R}{\mathbb C}$ associated to the eigenvalues $i$, $-i$ respectively, are involutive.
\end{definition}

All CR manifolds in this paper will be smooth, closed, co-oriented and oriented, and all maps and 
tensors will also be  smooth  unless otherwise stated.

The Levi form of a CR manifold $(M,D,J)$ is the symmetric
bilinear tensor given by
\begin{eqnarray*}
D\times D& \longrightarrow & TM/D\\
(u,v)&\longmapsto & [U,JV]/\sim\;\;,
\end{eqnarray*}
where $U, V$ are local sections of $D$ extending $u,v\in T_xM$, and we consider the class of 
the above Lie bracket at $x$ in the quotient real line bundle $TM/D$ (which is oriented, so we can make sense of
positive and negative values).
The Levi form keeps track of the behavior of the distribution $D$.
 Its vanishing  is equivalent to $D$ integrating into a foliation $\mathcal{F}$,
  in which case we speak of a Levi-flat CR manifold. The other extreme case
   is that of strictly pseudo-convex (resp. pseudo-concave)  CR manifolds,
for which the  Levi form is strictly positive (resp. negative); in particular the distribution 
$D$ of such CR manifolds is a contact
     distribution.

Let $(M,D,J)$ be a CR manifold and $(M',D',J')$ either
 a CR manifold or a complex manifold (in which case $D'=TM'$). A map $\phi\colon M\rightarrow M'$
 is CR  if $\phi_*D\subset D'$ and $\phi_*\circ J=J'\circ \phi_*$. A  CR vector bundle is a complex
vector bundle $\pi\colon E\rightarrow (M,D,J)$ defined by CR transition maps.

\subsection{Positivity of CR line bundles and CR sections.}

Our goal is finding CR line bundles $L\rightarrow  M$ with plenty of CR sections, so among them we
 have linear systems of CR sections whose leafwise holomorphic jets solve appropriate transversality problems.
 These would be our generic CR linear systems. A natural condition to impose is that of positivity of $L$ along $D$; in the
Levi-flat case and if all leaves are compact -so from the
differential viewpoint the foliated manifold $(M,{\mathcal{F}})$ is a
mapping torus- $L$ restricts to each leaf to a positive line bundle, and therefore large enough powers
 admit generic holomorphic linear systems. One expects suitable $S^1$-families of such generic holomorphic
 linear systems to fit into generic CR linear systems of  $L^{\otimes k}$.

Among CR linear systems generic ones are open, since they are defined by transversality conditions.
 Roughly speaking the way in which one would think of producing them is by breaking the problem into two parts: firstly,
 finding linear systems 
which are already generic and very close to be CR (genericity can also be extended to not necessarily
 CR linear systems). Secondly, projecting the previous  generic nearly CR linear system into a CR linear system by 
solving the corresponding
 tangential Cauchy-Riemann equation, and making sure that the correction is small enough so that genericity is preserved.
 The first part of the problem can be solved: the assumption of positivity along $D$ implies that the curvature
 of a compatible Hermitian connection on $L$ divided by $-2\pi i$, makes $(M,D)$ into a 2-calibrated
 manifold for which $J$ is a compatible almost complex structure. The theory developed in \cite{Ma04}
 grants the existence of generic nearly CR sections.  What cannot be granted in general  is the
 existence of a solution to the tangential Cauchy-Riemann problem close enough to the nearly CR section;
 for Levi-Flat manifolds  this should be possible by using results of Ohsawa-Sibony \cite{OS00}.
 We prefer to choose another approach in which we assume that $M$ has a CR embedding into a (compact)
 complex manifold $X$, and $L$ is the restriction of a positive holomorphic line bundle $\mathcal{L}\rightarrow X$.
 Moreover, rather than looking at arbitrary CR sections of $L^{\otimes k}$, we will be only
 interested in those which are restriction of holomorphic sections of ${\mathcal{L}}^{\otimes k}$,
 k $\gg$ 1, the advantage being that for the latter line bundles we have a lot of control for the projection 
of a section onto its holomorphic part. Still, one needs to show that among these very particular CR sections 
genericity can be achieved. In the language of \cite{Ma04} what we explained is why we do not 
know how to adapt intrinsic approximately holomorphic theory to the CR setting, and we need to
 impose the existence of ${\mathcal{L}}\rightarrow X$ and explore  whether relative 
approximately holomorphic theory can be adapted to $(L,M)\hookrightarrow ({\mathcal{L}},X)$.

For simplicity we will assume $X$ to be some projective space and $\mathcal{L}$ to be the hyperplane
 line bundle ${\mathcal{O}}(1)$.

\subsection{Projective CR manifolds.} \begin{definition}\label{def:proj} A CR manifold $(M,D,J)$ is projective
 if it admits a CR embedding into some ${\mathbb C}{\mathbb P}^N$.
It is called embeddable if it admits a CR embedding into some
${\mathbb C}^N$.
\end{definition}

A large number of projective CR manifolds are provided by the following embedding theorems.

\begin{theorem}\label{thm:embed}(Boutet de Monvel \cite{BdM75}) Any (compact and oriented)
 strictly pseudo-convex (resp. pseudo-concave) CR manifold of
 dimension bigger or equal  than five is embeddable.
\end{theorem}

A recent result of Marinescu and Yeganefar \cite{MY04}  states that any Sasakian
manifold is embeddable. Sasakian manifolds are strictly pseudo-convex. 3-dimensional ones
 -to which theorem \ref{thm:embed} does not apply-
   are also covered by the aforementioned result of Marinescu and Yeganefar.

\begin{theorem} \label{thm:embed2} (Ohsawa, Sibony \cite{OS00})\label{thm:OS} Let $(M,D,J)$
 be a (compact) Levi-flat CR manifold and $L\rightarrow M$ a positive line bundle.
  Let h be any natural number. Then there exists another natural number d(h)
    and  a CR embedding  $j\colon M\hookrightarrow {\mathbb C}{\mathbb P}^{d(h)}$ of class $C^h$.
\end{theorem}

\subsection{1-jets and Lefschetz pencil genericity condition.}
Let $(M,D,J)$ be a CR manifold (not necessarily  projective). Let $(M'J')$ be either
 a complex or a CR manifold, and let us denote in either case by $D'$ the
  maximal complex distribution of $TM'$.

\begin{definition}\label{def:1jets}
The bundle of CR $1$-jets of CR maps from $M$ to $M'$, denoted by
${\mathcal{J}}^1_{\rm{CR}}(M,M')\rightarrow M$, is defined to  be
the bundle over $M$ whose fiber over $x$ is
\[{\mathcal{J}}^1_{\rm{CR}}(M,M')_x:=\{(y, h)\;|\;y\in M', h \in {\rm{Hom}}_{\mathbb C}(D_x,D'_y)\}.\]

Inside this bundle there is a distinguished submanifold

\[\Sigma=\{(x,y,h)\in {\mathcal{J}}^1_{\rm{CR}}(M,M')\;|\;h=0\}.\]

The CR 1-jet of a CR map $\phi\colon M\rightarrow M'$ is by definition
 $j^1_{\rm{CR}}\phi :=(\phi,d_D\phi)$, where $d_D\phi$ is the  restriction
  of the differential $d \phi$ to $D$. It is a section of ${\mathcal{J}}^1_{\rm{CR}}(M,M')$.
\end{definition}

If $(M,D,J)$ is a projective CR manifold
 there is a second relevant complex distribution associated to $M$: the complex envelope
 of its tangent bundle, denoted by $TM^J$, where at each point $x\in M$ we put the smallest complex
  subspace of $T_x{\mathbb C}{\mathbb P}^N$  containing $T_xM$ . For some $\epsilon>0$ small enough, we use
 the Fubini-Study metric to extend both $D$ and $TM^J$ by parallel
 transport along normal geodesics
  to complex distributions defined in ${\mathcal{N}}_\epsilon(M)$, the neighborhood of $M$ of radius $\epsilon$.
   The extended distributions will be also denoted  by $D$ and $TM^J$ whenever there is no risk of confusion.

 We have the corresponding  bundle ${\mathcal{J}}^1_{TM^J}({\mathcal{N}}_\epsilon(M),M')\rightarrow {\mathcal{N}}_\epsilon(M)$
  of holomorphic 1-jets along $TM^J$ of holomorphic maps from  ${\mathcal{N}}_\epsilon(M)$ to $M'$.
   We denote  by $\Sigma^J$ the submanifold corresponding to the holomorphic 1-jets along
    $TM^J$ whose degree 1 homogeneous component is vanishing. Similarly, the holomorphic 1-jet along $TM^J$ 
of a holomorphic map
 $\Phi\colon {\mathcal{N}}_\epsilon(M)\rightarrow M'$ is by definition
 $j^1_{TM^J}\Phi :=(\Phi,d_{TM^J}\Phi)$. It is a section of ${\mathcal{J}}^1_{TM^J}({\mathcal{N}}_\epsilon(M),M')$.

\begin{definition}\label{def:Lpencilstr} Let $j\colon (M,D,J)\hookrightarrow {\mathbb C}{\mathbb P}^N$ be
a projective CR manifold.  A (CR) Lefschetz pencil structure (of degree k) on $M$ is defined to be
 a  pencil of  (degree k) hypersurfaces of ${\mathbb C}{\mathbb P}^N$, so that
  any two hypersurfaces intersect transversely in the base locus $\tilde{B}$, and the following conditions hold:

\begin{enumerate}
\item  $\tilde{B}$ is transverse to $M$ and therefore  $B:=\tilde{B}\cap M$
 is a real codimension four CR submanifold of $M$.
\item Let $\Phi\colon {\mathbb C}{\mathbb P}^N\backslash \tilde{B}\rightarrow {\mathbb C}{\mathbb P}^1$ denote
 the holomorphic map associated to the pencil and $\phi\colon M\backslash B\rightarrow{\mathbb C}{\mathbb P}^1$
  its restriction to  $M\backslash B$. Then  
\[j^1_{\rm{CR}}\phi\colon  M\backslash B\rightarrow  {\mathcal{J}}^1_{\rm{CR}}(M\backslash B,{\mathbb C}{\mathbb P}^1)\]
 is transverse along $D$ to $\Sigma$.
\item  The holomorphic 1-jet along $TM^J$
   \[j^1_{TM^J}\Phi\colon {\mathcal{N}}_\epsilon(M) \backslash \tilde{B}\rightarrow
 {\mathcal{J}}^1_{TM^J}({\mathcal{N}}_\epsilon(M)\backslash \tilde{B},{\mathbb C}{\mathbb P}^1)\]
 is  transverse along $M$ to $\Sigma^J$.
\end{enumerate}
We will also call the pencil of hypersurfaces giving a Lefschetz
pencil structure a
  generic rank one (CR) linear system of the CR line bundle ${\mathcal{O}(k)}_{\mid M}:=j^*{\mathcal{O}}(k)$.
\end{definition}

We denote the Lefschetz pencil structure by the triple
$(\phi,B,\Delta)$, where $\Delta$ are the points of $M$ where the
derivative $d_D\phi$ vanishes. Any point $x\in \Delta$ is called a singular point of the pencil.
 Points in $B$ are base points, and points in
 $M\backslash (B\cup \Delta)$ are called regular. Similarly,
 $a\in {\mathbb C}{\mathbb P}^1\backslash \phi(\Delta)$ is called a regular value
  and singular otherwise.
Notice that at a singular point of the pencil the map $\phi$ fails to be a submersion.

\subsubsection{Transversality along distributions.} Recall that transversality along $D$ to $\Sigma$ of a section
$\sigma\colon M
\rightarrow{\mathcal{J}}^1_{\rm{CR}}(M,{\mathbb C}{\mathbb P}^1)$ is
defined as follows \cite{Ma04}:
 the pullback of $D$ by the projection  ${\mathcal{J}}^1_{\rm{CR}}(M,{\mathbb C}{\mathbb P}^1)\rightarrow M$
   defines a distribution $\hat{D}$ in ${\mathcal{J}}^1_{\rm{CR}}(M,{\mathbb C}{\mathbb P}^1)$.
   At every point $x$ such that $\sigma(x)\in \Sigma$  one requires
    \[ \hat{D}_{\sigma(x)}=T_{\sigma(x)}\Sigma\cap \hat{D}_{\sigma(x)} + T_x\sigma\cap \hat{D}_{\sigma(x)},\]
 where $T_x\sigma$ is the tangent space of the graph of $\sigma $ at $\sigma(x)$.

The above definition extends in the obvious way to transversality along $D$ to any submanifold
 $S$ of a bundle $E\rightarrow M$.

Transversality along $M\subset {\mathcal{N}}_\epsilon(M)$ is defined at the points of $M$ as we did for 
transversality along $D$, but with $TM\subset T{\mathcal{N}}_\epsilon(M)_{\mid M}$ playing 
the role of $D$. Since $TM\subset {T{\mathcal{N}}_\epsilon(M)}_{\mid M}$ is integrable, 
transversality along $M$ of $j^1_{TM^J}\Phi$
  is equivalent to  usual transversality for the restriction
  \[j^1_{TM^J}\Phi\colon  M\backslash B \rightarrow
   j^*{\mathcal{J}}^1_{TM^J}({\mathcal{N}}_\epsilon(M)\backslash
   \tilde{B},{\mathbb C}{\mathbb P}^1).\]

Transversality along distributions is stronger than usual transversality.

From the transversality assumption along $D$ on
$j^1_{\rm{CR}}\phi$  we conclude that $\Delta$ is a 1-dimensional submanifold transverse to $D$.
 A dimension count implies that transversality of $j^1_{TM^J}\Phi$ along $M$ to $\Sigma^J$
  is equivalent  to the section not intersecting  $\Sigma^J$.
   Therefore in the points of $\Delta$ the derivative $d\phi\colon TM\rightarrow T{\mathbb C}{\mathbb P}^1$
   cannot vanish, because being $\Phi$ holomorphic  that would imply the vanishing of $d_{TM^J}\Phi$.
    Thus $\phi(\Delta)\subset {\mathbb C}{\mathbb P}^1$ is an immersed curve, which
    is why condition 3 in definition \ref{def:Lpencilstr} is required.

\begin{remark}\label{rem:class} If the CR embedding $j\colon M\hookrightarrow {\mathbb C}{\mathbb P}^N$ 
is of class $C^h$, h$\geq$ 2,
  definition \ref{def:Lpencilstr} also makes sense because transversality (along any distribution) is a $C^1$-notion.
\end{remark}

Our first main result is the following:

\begin{theorem}\label{thm:pencil} Let $(M^{2n+1},D,J)$ be a
  projective CR manifold of class $C^h$, h $\geq$ 3.
     Then for any  k large enough  $M$ admits degree k Lefschetz pencil structures.
       Let $W$ be the fiber of $\phi$ over any regular value,
        compactified by adding $B$. Then the morphisms at the level
         of homotopy (and homology groups)
         $l_*\colon \pi_i(W)\rightarrow \pi_i(M)$ induced by the inclusion are
          isomorphisms for 0 $\leq$ i $\leq$ n-2, and an epimorphism for i $=$ n-1.
\end{theorem}

Theorem \ref{thm:pencil} is a result about the existence of generic rank one
 linear systems of the bundles ${\mathcal{O}(k)}_{\mid M}\rightarrow M$, for k $\gg$ 1.

It is not possible in general to give a normal form for $\phi$ in
neighborhoods of points of  $\Delta$  and $B$. In the Levi-flat setting 
normal forms exist.

\begin{proposition}\label{pro:nforms}(see also  \cite{IM04},\cite{Ma05}) Let $(\phi,B,\Delta)$ be a
 Lefschetz pencil structure on a projective Levi-flat CR manifold. Then we have the following:
\begin{enumerate}
\item For any $x\in B$  there exist CR coordinates $z^1,\dots,z^n,s$ centered
 at $x$ and  a holomorphic chart of ${\mathbb C}{\mathbb P}^1$, such that
 \[B\equiv z^1=z^2=0\] and \[\phi(z^1,\dots,z^n,s)=z^2/z^1.\]
\item For any $x\in \Delta$ there exists CR coordinates $z^1,\dots,z^n,s$
centered at $x$, a holomorphic chart of ${\mathbb C}{\mathbb P}^1$ and
a real map $t$ of class $C^{h-1}$ such that
\[\phi(z^1,\dots,z^n,s)={(z^1)}^2+\cdots +{(z^n)}^2+t(s),\; t(0)=0,\;
t'(0)\neq 0.\]
\end{enumerate}
\end{proposition}

Due to the existence of normal forms it is possible to adapt results of \cite{Ma05} and know 
more about the intersection of a regular fiber of the pencil with a leaf of the foliation.

\begin{theorem}\label{thm:hyp} Let $(M,{\mathcal{F}},J)$ be a  projective Levi-flat CR
 manifold of dimension bigger or equal than five,
and let $(\phi,B,\Delta)$ be one of the pencils provided by theorem
\ref{thm:pencil}. Then for any regular fiber $W$ of $\phi$
(compactified adding the base points) and for every leaf $F$ of
$\mathcal{F}$, the intersection $F\cap W$ is connected.
As a consequence the inclusion $l\colon
(W,{\mathcal{F}}_W)\hookrightarrow (M,\mathcal{F})$,
 where ${\mathcal{F}}_W$ is the induced foliation, descends to a homeomorphism of leaf spaces.
\end{theorem}

Theorem \ref{thm:hyp} can also be understood as a leafwise Lefschetz hyperplane theorem
for the $\pi_0$, and for  sections (degree k hypersurfaces) which
fit into a Lefschetz pencil for  $(M,{\mathcal{F}},J)$. Notice as well how unexpected this result
 is for non-compact leaves (some of which can be dense!).

We would like to single out theorem \ref{thm:pencil} together with proposition \ref{pro:nforms} for
3-dimensional  projective Levi-flat CR manifolds. 

Let $M^3$ be a closed, orientable 3-manifold endowed with a smooth co-orientable foliation
  $\mathcal{F}$ by surfaces. We recall the following elementary result.

\begin{lemma} $(M^3,\mathcal{F})$ admits the structure of a Levi-flat CR manifold (with CR foliation $\mathcal{F}$)
with a positive CR line bundle, if and only if $\mathcal{F}$ is a taut
foliation.
\end{lemma}
\begin{proof}  The existence of a Levi-flat CR structure with a positive CR Hermitian line bundle
 $(L,\nabla)$ with compatible connection clearly implies tautness. 
Following Sullivan \cite{Su76} tautness is equivalent to the existence
 of a closed 2-form restricting to a leafwise area form.
 In our case the curvature of $\nabla$ divided by $-2\pi i$  has this property.

Conversely, if the foliation is taut we do have an integral closed
2-form $\omega$ which is non-degenerate when restricted to each
leaf.
 Let $J$ be a leafwise compatible almost complex structure. Then $J$ is integrable, for the
 leaves are 2-dimensional, and $\omega$ is of type (1,1) with respect to $J$. Therefore  we can construct
 $(L,\nabla)$, the corresponding pre-quantum line bundle, which is a Hermitian line bundle with
  connection whose curvature satisfies $F_\nabla=-2\pi i\omega$. Being $\omega$ of type (1,1) the bundle
  admits a leafwise holomorphic structure, and hence a CR one. Positivity is also clear.
\end{proof}

\begin{corollary}\ Let $(M^3,\mathcal{F})$ be a closed, orientable 3-manifold endowed with a
co-orientable taut foliation, and fix $h\in \mathbb{N}$, $h\geq 3$. Then there exists a CR structure in
$(M^3,\mathcal{F})$  and maps $\phi\colon M^3\rightarrow
{\mathbb C}{\mathbb P}^1$ with the following properties:
\begin{enumerate}
\item  $\phi$ is of class $C^h$ and leafwise holomorphic.
\item The restriction of $\phi$ to each leaf is a branched cover with index 2 singular points.
\item The leafwise singular sets fit into a transverse link $\Delta$ of class $C^{h-1}$.
\item Around each point $a\in \Delta$ there exist local CR coordinates $(z,s)$ and a complex coordinate
 in ${\mathbb C}{\mathbb P}^1$, so that  $\phi(z,s)=z^2+t(s)$, where $t$ is $C^{h-1}$ and $t(0)=0$, $t'(0)\neq 0$.
\end{enumerate}
\end{corollary}

\begin{corollary}  Let $(M^3,{\mathcal{F}})$ be a closed and orientable 3-manifold endowed with a
co-orientable foliation. Then ${\mathcal{F}}$ is taut is and only if
there exist a leafwise complex structure for which
$(M^3,{\mathcal{F}},J)$ admits a (CR) Lefschetz pencil structure.
\end{corollary}

Also observe that by applying theorem \ref{thm:hyp} inductively we conclude that any
 projective Levi-flat CR manifold contains a taut 3-dimensional foliation whose leaves
 are complex, and such that the inclusion descends to a homeomorphism of leaf spaces. 
Thus leaf spaces of projective Levi-flat CR manifolds are no more complicated than those of 3-dimensional taut foliations. 

\subsection{Higher order jets.}
Regarding higher order genericity, note that if $(M,D,J)$ is a CR manifold, $M'$ either a CR or complex manifold  and
 $\phi\colon M\rightarrow M'$ a CR map, it is not possible in general to define the second CR jet of $\phi$;
 there is no intrinsic way
of taking a second derivative along a general distribution $D$. Therefore, to consider higher order jets we must assume that
 $(M,{\mathcal{F}},J)$ is Levi-flat. 

 On the bundle of foliated holomorphic r-jets of CR maps from $(M,{\mathcal{F}},J)$ to $M'$ one has 
the  CR Thom-Boardman stratification $\mathcal{T}^{\mathcal{F}}$, which is the union of the usual holomorphic Thom-Boardman
stratifications \cite{Bo67} over each leaf $F\in \mathcal{F}$. A CR map  $\phi$  is defined to be r-generic if its
foliated holomorphic r-jet (i.e. its CR r-jet) is leafwise  transverse to
$\mathcal{T}^{\mathcal{F}}$. For such a map the pullback of each
stratum $\Sigma^{\mathcal{F}}_I(\phi)$ would be transverse to the
leaves; its intersection with each leaf $F$  would be the
holomorphic Thom-Boardman stratum $\Sigma_I(\phi_{\mid F})$.

 Let us for simplicity  forget about holomorphic functions and consider
  the foliated/leafwise genericity problem in the smooth setting.
  A strategy to solve transversality problems for foliated smooth jets
   is to use the canonical submersion from the bundle of r-jets to the bundle
   of foliated r-jets \cite{Be02},  to pull back the  leafwise Thom-Boardman
   stratification to the bundle of r-jets. Leafwise r-genericity is equivalent
    to transversality along the leaves of $\mathcal{F}$ of the (full) r-jet to the pulled back stratification.

To develop a similar strategy to solve  the leafwise holomorphic genericity
 problem, we would need to embed $M$ inside a complex manifold $X$, and transfer
  the problems for leafwise holomorphic
  jets to problems for full holomorphic jets. The next step would be to
  solve the foliated (strong) transversality problem in the bundle of full
   holomorphic r-jets. This is always possible locally, but there is no
    reason why a global solution should exist.

If $(M,{\mathcal{F}},J)$ is Levi-flat and possess a positive CR line
bundle, theorem
 \ref{thm:embed2} implies that for any natural number $h$  we obtain a
  leafwise holomorphic embedding in ${\mathbb C}{\mathbb P}^{d(h)}$ of class $C^h$. This
  projective space plays the role of the ambient complex space $X$.

The target space of the r-generic maps we will construct  will not
be an arbitrary complex manifold,
 but complex projective space, for our maps will come again from linear systems.

\begin{theorem}\label{thm:gener} Let $(M^{2n+1},{\mathcal{F}},J)$ be a Levi-flat
 CR manifold endowed with a positive CR line bundle $L\rightarrow M$. Fix
 h,r $\in {\mathbb N}$, r $\leq$ h-2. Then for any integer m  we can find the following: 
\begin{enumerate}
\item A Levi-flat CR submanifold $B$ of real codimension 2m+2 and class $C^h$.
\item A CR map $\phi\colon M\backslash B \rightarrow {\mathbb C}{\mathbb P}^m$ -which is the restriction
of a holomorphic map
 $\Phi\colon {\mathbb C}{\mathbb P}^N\backslash {\mathcal{B}}\rightarrow {\mathbb C}{\mathbb P}^m$- such that
  its CR r-jet  is leafwise transverse to the Thom-Boardman stratification of
    the bundle of CR r-jets (which is  by Levi-flat strata of class $C^{r-h}$).
   Therefore the pullback of each stratum is a
   Levi-flat submanifold of the expected codimension and of class $C^{r-h}$.
\item The 1-jet of $\Phi$ along $TM^J$ is transverse along $M$ to
the stratification of
${\mathcal{J}}_{TM^J}^1({\mathcal{N}}_\epsilon(M),{\mathbb C}{\mathbb P}^m)$,
whose strata are defined according to the rank of the degree 1 homogeneous
component of $\sigma\in
{\mathcal{J}}_{TM^J}^1({\mathcal{N}}_\epsilon(M),{\mathbb C}{\mathbb P}^m)$.
\end{enumerate}
\end{theorem}

Observe that the outcome of theorem \ref{thm:gener} is a generic CR map $\phi$ defined away from a CR 
submanifold (properties (1) and (2)), together with a requirement on the full
derivative of $\Phi$  analogous to the third one in the
definition of a Lefschetz pencil structure.  The latter condition gives a
lower bound  -depending on the dimensions m,n- on the rank
$d\phi$, this providing  supplementary information about the
differential of $\phi$ in the direction transverse to $D$ in those
points where the rank along $D$ drops enough.

Notice as well that theorem \ref{thm:gener0} in the introduction is a simplified version of theorem \ref{thm:gener}.

\section{Dual geometry of projective CR manifolds}\label{sec:sec3}

For any projective manifold $X\hookrightarrow
{\mathbb C}{\mathbb P}^N$  generic linear systems of ${\mathcal{O}}(1)$
 can be constructed  by applying basic transversality results to
 its dual variety
  $X^*\hookrightarrow {\mathbb C}{\mathbb P}^{N*}$. It might happen that the dual variety is
  not a divisor,
   but it is indeed a divisor if we twist the embedding (compose with
   the k-th Veronese embedding, for any k $\geq$ 2). So in the worst case
    there exist always degree 2 Lefschetz pencils  \cite{SGA7II}.

Let $(M,D,J)\hookrightarrow {\mathbb C}{\mathbb P}^N$ be a projective
CR variety. We define its dual set
\[ M^*=\{H\in {\mathbb C}{\mathbb P}^{N*}\;|\;D_x\subset H_x \;{\rm{for}}\,\, {\rm{some}}\; x\in M\}.\]
By pulling back the usual diagram of flag varieties, it can be shown that there exists
 $\pi\colon S_D\rightarrow M$ a smooth (resp. $C^{h-1}$ for $C^h$-embeddings
  of Levi-flat manifolds) fiber bundle of real dimension 2N-1, and a smooth
   (resp. $C^{h-1}$) map $\nu\colon S_D\rightarrow  {\mathbb C}{\mathbb P}^{N*}$ -the dual map- such  that $M^*=\nu(S_D)$.

The difference with respect to the complex setting -where $X^*$ is known
 to be a (singular) complex variety because it is the image of
  a complex manifold by a complex map-  stems  from the fact that $S_D$
   has only the structure of smooth manifold,
   and $\nu$ is just a smooth map. When $(M,{\mathcal{F}},J)$ is Levi-flat the situation is slightly
better, but still not good enough to do geometry with the dual
variety. Indeed, the fiber bundle $S_D$ is a Levi-flat manifold and
the dual map
 $\nu\colon S_D\rightarrow  {\mathbb C}{\mathbb P}^{N*}$ is a CR map. If $y\in S_D$
 is a regular point for the leafwise dual map, then it is easy to check that it
  is a regular point for the full dual map. Hence, $\nu$ fails to be regular
  where the leafwise dual map is degenerate. That might lead to some control
   on the singular points $S^*$ of $M^*$, which in turn might imply
    the existence of pencils of hyperplanes
   avoiding $S^*$ and intersecting the regular part transversely, i.e. pencils fulfilling 
conditions (1) and (2) in definition \ref{def:Lpencilstr}.  If $(M,{\mathcal{F}},J)$ is a Levi-flat manifold of class at least $C^2$
      it cannot have codimension one in ${\mathbb C}{\mathbb P}^{N}$, $N\geq$ (or Lipschitz if $N\geq 3$) \cite{CSW04,CSW05}. Then
       inside of $S_D$ there is a real codimension $2$ sub-bundle
       $S_0$ corresponding to those hyperplanes, that as well as
       containing some $D_x$, also contain $T_xM^J$. Let $S_0^*:=\nu(S_0)$.
        In order to construct degree 1 Lefschetz pencil structures
          as defined in \ref{def:Lpencilstr} we would further need to ask the pencils 
          fulfilling conditions (1) and (2), to miss $S^*_0$.  The difficulty
             comes from the fact that it is not clear that $S_0$
             is a Levi-flat manifold (it is not clear that the leaves
              of the obvious codimension one foliation are complex), so once more we 
cannot say much about its image $S^*_0$. In particular we cannot argue that the set of pencils not
              intersecting $S^*_0$ is non-empty.

We can just use Sard's theorem to
 argue that $M^*$ has measure zero and hence that hyperplane sections
  do exist. We recall that the topology of the hyperplane section
  can only be related to that of $M$  when the  Levi form has some
  degeneracy, in which case the results of Ni and Wolfson apply \cite{NW03}.

One might think heuristically  of theorems \ref{thm:pencil} and
\ref{thm:gener}
 as a manifestation of the existence some sort of
dual geometry  for the re-embeddings provided by the k-th Veronese
maps for k very large.  That is,  the corresponding dual sets
$M^*_k$ are in ``most of its points'' close to be stratified
varieties; similarly whenever $M$ has real codimension bigger than
one, and for k $\gg$ 1, the image of  $S_{0,k}$ should be thought of
being close to be a stratified variety with complex strata of
complex codimension at least two. Therefore, a generic pencil of hyperplanes
should be able to avoid both the closed strata of $M^*_k$ and the
image of $S_{0,k}$. However, we were not able to come up with a
proof of this fact and hence with a geometric proof of theorems
\ref{thm:pencil} and \ref{thm:gener}.

In a similar vein the results about the
 topology of the smooth fiber of a Lefschetz pencil should not come as a surprise.
  They coincide with the results of \cite{NW03} for Levi-flat manifolds.
   And one should bear in mind that approximately holomorphic theory is based on the study of the
   CR manifold at a very small scale, where it looks like a Levi-flat one.

As we just discussed dual geometry fails in general to produce Lefschetz pencils  for 
projective CR manifolds, due to the lack of structure on the dual set $M^*$. In \cite{Ma08}
 we studied real hypersurfaces $M$ of projective space for which the dual map 
 $\nu\colon S_D\cong M\rightarrow {\mathbb C}{\mathbb P}^{N*}$ is an immersion.
 It is proved that injectivity of the differential of $\nu$ at any $x\in M$ is equivalent
 to the real shape operator restricted to $D_x=T_xM\cap JT_xM$ being non-degenerate. We also
 showed that for such hypersurfaces a pencil of hyperplanes $\mathbb{L}\subset {\mathbb C}{\mathbb P}^{N*}$ 
was a Lefschetz pencil as in definition \ref{def:Lpencilstr} (for real hypersurfaces condition 
 3 is void), if and only if $\mathbb{L}$ was transverse to $M^*$ (the later being a compact immersed 
hypersurface, transversality at each point is understood as transversality to all branches 
through the point). A dimension count and Sard's theorem easily gives the existence of complex 
lines transverse to $M^*$, therefore dual geometry provides Lefschetz pencils for these hypersurfaces.

A real hypersurface $M\hookrightarrow  {\mathbb C}{\mathbb P}^{N}$ for which the shape operator
 restricted to the CR distribution $D$ is strictly positive at each point, is called an strictly
 $\mathbb C$-convex hypersurface. These hypersurfaces are embeddable, since they miss hyperplanes. 
In particular one can always find pencils of hyperplanes whose base misses $M$, and among them 
Lefschetz ones. The corresponding CR function $\phi$ does not have base points and it is 1-generic,
 so we can talk of a CR Morse function. For strictly convex real hypersurfaces (inside Euclidean
 space) pencils of real hyperplanes whose base misses the hypersurface provide the simplest
 Morse function, i.e. one with two critical points, a maximum and a minimum; in particular
 they give a homeomorphism from the hypersurface to the sphere. Our CR Morse functions on 
strictly $\mathbb C$-convex hypersurfaces are analogs of the former: the singular set $\Delta$ 
is just a copy of $S^1$, and it is mapped diffeomorphically by $\phi$ to the boundary of the 
image of $\phi$, which is diffeomorphic to a closed disk (in the real case the image is a closed 
interval bounded by the image of the two critical points). It is also possible to reconstruct $M$
 from the CR Morse function to conclude that it is diffeomorphic to $S^{2N-1}$ \cite{Ma08} (a result already
 known but just at the topological level \cite{Ho}).

\section{Approximately holomorphic geometry for projective CR manifolds}\label{sec:sec4}

Let us fix some notation. Given $s_k\colon
{\mathbb C}{\mathbb P}^N\rightarrow
\underline{\mathbb C}^{m+1}\otimes {\mathcal{O}}(k)$, the restriction  of $s_k$  to $M$
  will be denoted  by  $\tau_k$, and its zero set is the base locus $B$. The projectivization
  of $\tau_k$ to $M\backslash B$ will
be denoted by $\phi_k$.

   The  holomorphic vector bundles  $\underline{\mathbb C}^{m+1}\otimes {\mathcal{O}}(k)$
    -which will be also denoted by $E_k$- carry a natural connection $\nabla_k$
     coming from the flat one in $\underline{\mathbb C}^{m+1}$ and the connection
      in ${\mathcal{O}}(k)$ associated to the Fubini-Study form. We will use the same
       notation for the restriction of $\nabla_k$ to ${E_k}_{\mid M}$ if there
        is no risk of confusion.

For a  projective Levi-flat CR manifold of class $C^h$, we let
${\mathcal{J}}_{\rm{CR}}^r(M,{\mathbb C}{\mathbb P}^m)$ denote the
bundle of CR $r$-jets of CR maps to
${\mathbb C}{\mathbb P}^m$. This is a bundle of class $C^{h-r}$ and it
inherits an obvious CR-structure. There is a leafwise Thom-Boardman
stratification $\mathbb P{\mathcal{T}}^{\mathcal{F}}$ whose strata are
Levi-flat CR submanifolds of class $C^{h-r}$. For any CR map $\phi$
to ${\mathbb C}{\mathbb P}^m$ of class $C^h$, its  CR r-jet
prolongation $j^r_{\rm{CR}}\phi$ is a CR section of class
$C^{h-r}$.

In order to prove points (1) and (2) of theorem \ref{thm:gener}, we need
to find suitable sequences of sections $s_k$ of $E_k$ such that (1)
the base locus of $\tau_k$ is a CR submanifold of the expected
dimension, and (2) $j^r_{\rm{CR}}\phi$ is transverse along
$\mathcal{F}$ to $\mathbb P{\mathcal{T}}^{\mathcal{F}}$. This two conditions -as we will see- can 
actually be stated in a more compact manner 
as the solution to a single transversality problem for $\tau_k$. What we will do is breaking the problem into three parts:
\begin{enumerate}
\item[(A)]  The aforementioned transversality problem can be linearized, i.e. 
the bundle ${\mathcal{J}}_{\rm{CR}}^r(M,{\mathbb C}{\mathbb P}^m)$, the Thom-Boardman stratification and
 the notion of CR jet. 
\item[(B)] The linearized problem has a CR solution $\tau_k$ -with projectivization $\phi_k$- 
which will be provided by a suitable version 
of approximately holomorphic geometry.
\item[(C)] The CR solution of the linearized problem also solves the original problem.
\end{enumerate}

We need to recall a number of notions and results from approximately holomorphic geometry.
Let  $s_k\in\Gamma(E_k)$.  Using $J$ the complex structure of ${\mathbb C}{\mathbb P}^N$ we can write

  \[\nabla s_k=\partial s_k+\bar{\partial}s_k,\; \partial s_k\in \Gamma(T^{*1,0}{\mathbb C}{\mathbb P}^N\otimes E_k),\;
\bar\partial s_k\in \Gamma(T^{*0,1}{\mathbb C}{\mathbb P}^N\otimes
E_k).\]

Similarly, given $\tau_k\in \Gamma({E_k}_{\mid M})$
 the restriction of $\nabla\tau_k$ to $D$ can be written
  \[\nabla_D\tau_k=\partial\tau_k+\bar{\partial}\tau_k,\; \partial\tau_k\in \Gamma(D^{*1,0}\otimes {E_k}_{\mid M}),\;
\bar\partial\tau_k\in \Gamma(D^{*0,1}\otimes {E_k}_{\mid M}).\]

Let $g$ denote the Fubini-Study metric and let $g_k$  denote
 the rescaled metric $kg$. We use the same notation for the restriction of these metrics to $M$.

\begin{definition}\label{def:sucah}  A sequence of sections   $s_k$ of $E_k$ is
approximately $J$-holomorphic (or approximately holomorphic or
simply A.H.), if positive constants ${(C_j)}_{j\geq 0}$ exist such
that for all $k\gg 1$

\[  |\nabla^j s_k|_{g_k}\leq  C_j,\;\; |\nabla^{j-1} \bar{\partial} s_k |_{g_k} \leq C_j k^{-1/2}. \]
If in an A.H. sequence the sections $s_k$ are
holomorphic we speak of a uniformly bounded  sequence of holomorphic sections.

Similarly, a sequence of sections   $\tau_k$ of ${E_k}_{\mid M}$ is
approximately $J$-holomorphic (or approximately holomorphic or
simply A.H.), if positive constants ${(C_j)}_{j\geq 0}$ exist such
that for all $k\gg 1$

\[  |\nabla^j\tau_k|_{g_k}\leq  C_j,\;\; |\nabla^{j-1} \bar{\partial}\tau_k |_{g_k} \leq C_j k^{-1/2}. \]
If the sections in the sequence are CR  we say that it is a uniformly
bounded  sequence of CR sections.
\end{definition}

The previous definition can also be given requiring control on a finite number of
 covariant derivatives, so we have $C^h$-A.H. sequences of sections when inequalities hold for 
 $j=0,\dots, h$.

\subsection{Linearization of the bundles of CR jets, the Thom-Boardman stratification, and the notion of CR r-jet}\label{ssec:lin1}

In this subsection we address point (A) in our strategy by presenting a linearized version of r-genericity 
for a sequence of sections of ${E_k}_{\mid M}$ which are not necessarily CR. 

Over the CR manifold $M$ we define the sequence of vector bundles of pseudo-holomorphic r-jets as
\begin{equation}\label{eqn:djets}
{\mathcal{J}}^r_D{E_k}_{\mid M}:=(\sum_{j=0}^r D^{*1,0}\odot\cdots^{(j)}\cdots\odot D^{*1,0})\otimes {E_k}_{\mid M},
\end{equation}
where $\odot$ denotes the symmetric product. These bundles carry a natural connection $\nabla_{k,r}$ and a metric (see
subsection 5.2 in \cite{Ma04} for more details).

The pseudo-holomorphic r-jet of $\tau_k$ will be a section of ${\mathcal{J}}^r_D{E_k}_{\mid M}$
 defined by induction: let   $j^{r-1}_D\tau_k \in {\mathcal{J}}^{r-1}_D{E_k}_{\mid M}$ be the
pseudo-holomorphic (r-1)-jet of $\tau_k$.
 It has homogeneous components of degrees $0,1,\dots,r-1$. We will denote
  the homogeneous component of degree $j\in \{0,\dots,r-1\}$ by
   $\partial^j_{\rm{sym}}\tau_k\in \Gamma(({D^{*1,0}})^{\odot j}\otimes {E_k}_{\mid M})$.
The connection $\nabla_{k,r-1}$ is actually a direct sum of connections defined
 on the direct summands $({D^{*1,0}})^{\odot j}\otimes {E_k}_{\mid M}$, $j=0,\dots,r-1$.
  For simplicity and if there is no risk of confusion, we will use the same
   notation for the restriction of $\nabla_{k,r-1}$ to each of the summands.
The restriction of
$\nabla_{k,r-1}\partial^{r-1}_{\rm{sym}}\tau_k$ to $D$ defines a
section \[\nabla_{k,r-1,D}\partial^{r-1}_{\rm{sym}}\tau_k\in
\Gamma(D^*\otimes({D^{*1,0}})^{\odot r-1}\otimes {E_k}_{\mid M}).\]
 For each $x\in M$ it is a form on $D$ with values in the complex vector space
 $({D^{*1,0}})^{\odot r-1}\otimes {E_k}_{\mid M}$. Therefore, we can consider
  its (1,0)-component \[\partial \partial^{r-1}_{\rm{sym}}\tau_k\in
\Gamma(D^{*1,0}\otimes({D^{*1,0}})^{\odot r-1}\otimes {E_k}_{\mid
M}).\] By applying the symmetrization map

\[{\rm{sym}}_j\colon {(D^{*1,0})}^{\otimes j}\rightarrow {(D^{*1,0})}^{\odot j}\]
we obtain $\partial^r_{\rm{sym}}\tau_k\in \Gamma(({D^{*1,0}})^{\odot r}\otimes {E_k}_{\mid M})$.

 \begin{definition}\label{def:psjets} Let  $\tau_k$ be a section of  $({E_k}_{\mid M},\nabla_k)$.
 The pseudo-holomorphic r-jet  $j^r_D\tau_k$ is a section of the bundle
 ${\mathcal{J}}^r_D{E_k}_{\mid M}$
  defined out of the pseudo-holomorphic (r-1)-jet by the formula
 \[j^r_D\tau_k:=(j^{r-1}_D\tau_k,\partial^r_{\rm{sym}}\tau_k).\]
 \end{definition}

Now we want to define the linearization of ${\mathcal{J}}_{\rm{CR}}^r(M,{\mathbb C}{\mathbb P}^m)$
 so that the ``projectivization'' of $j^r_D\tau_k$ -which is going to be the linearized r-jet of 
$\phi_k$- is a section of it. We will do it by gluing pieces defined very much as in equation \ref{eqn:djets}.
 
Let $Z^0,\dots,Z^m$ be the complex coordinates associated to the trivialization of
 $\underline{\mathbb C}^{m+1}$  and let $\pi\colon {\mathbb C}^{m+1}\backslash \{0\}\rightarrow {\mathbb C}{\mathbb P}^m$
 be the canonical projection. Consider the canonical affine charts
\begin{eqnarray}
 \nonumber\varphi_i^{-1}\colon U_i& \longrightarrow  & {\mathbb C}^m \\ \nonumber
 [Z_0:\dots:Z_m]& \longmapsto & \left(\frac{Z^1}{Z^0},\dots,\frac{Z^{i-1}}{Z^0},\frac{Z^{i+1}}{Z^0},\dots,\frac{Z^m}{Z^0}\right).
 \end{eqnarray}
For each chart $\varphi_i$, $i=0,\dots,m$, we consider the bundle

\begin{equation}\label{eqn:projjets}
{\mathcal{J}}_{D}^r(M,{\mathbb C}^m)_i:=(\sum_{j=0}^r {(D^{*1,0})}^{\odot j})\otimes \underline{\mathbb C}^m.
\end{equation}
On each bundle ${\mathcal{J}}_{D}^r(M,{\mathbb C}^m)_i$ we have a notion
 of pseudo-holomorphic r-jet as given in definition \ref{def:psjets}, where we use
 instead of ${E_k}_{\mid M}$ the trivial bundle $\underline{\mathbb C}^m$ with
  trivial connection and standard Hermitian metric associated to the frame
 $\xi_{i,1},\dots,\xi_{i,m}$ given by the above affine coordinates.

We gather a few key results concerning these vector bundles:

\begin{itemize}
\item By point (1) in proposition 6.1 in \cite{Ma04}, the vector bundles
${\mathcal{J}}_{D}^r(M,{\mathbb C}^m)_i$ can be glued to define the 
 fiber bundles ${\mathcal{J}}^r_D(M,{\mathbb C}{\mathbb P}^m)$  of pseudo-holomorphic
  r-jets of maps from $M$ to ${\mathbb C}{\mathbb P}^m$. One can put global metrics coming 
from $g_k$ in the base, though it is not necessary. The computations one needs to do are 
local so one can always assume that the sections belong to some ${\mathcal{J}}_{D}^r(M,{\mathbb C}^m)_i$,
 where we have metrics induced by $g_k$, the flat connection, and the standard
Hermitian metric.  
\item  According to point (2) in proposition 6.1  in \cite{Ma04}, 
given $\phi_k\colon M\rightarrow {\mathbb C}{\mathbb P}^m$  there exist a (unique) 
notion of pseudo-holomorphic r-jet extension 
\[j^r_D\phi_k\colon M\rightarrow {\mathcal{J}}^r_D(M,{\mathbb C}{\mathbb P}^m)\]
which is compatible with the notion of pseudo-holomorphic r-jet of definition  \ref{def:psjets} 
for the sections $\varphi_i^{-1}\circ\phi_k\colon M\rightarrow {\mathbb C}^m$. This is our linearized version of CR r-jet.
Note that neither $D$ nor $J$ have to be integrable. If $J$ is integrable the linearized notion
 of CR r-jet does not require $\phi_k$ to be a CR map.

\item   Define  ${\mathcal{J}}^r_D{E_k}_{\mid M}^*:={\mathcal{J}}^r_D{E_k}_{\mid M}\backslash Z_k$, where
  $Z_k$ denotes the sequence of strata of  ${\mathcal{J}}^r_D{E_k}_{\mid M}$ of r-jets whose
   degree 0 component vanishes. Then according to proposition 6.2   in \cite{Ma04},
$\pi\colon {\mathbb C}^{m+1}\backslash \{0\}\rightarrow {\mathbb C}{\mathbb P}^m$ induces bundle
 maps  $j^r\pi\colon {\mathcal{J}}^r_D{E_k}_{\mid M}^*\rightarrow {\mathcal{J}}^r_D(M,{\mathbb C}{\mathbb P}^m)$
such that  for any section $\tau_k$ of ${E_k}_{\mid M}$, in the points where it does not
 vanish and its projectivization $\phi_k$ is defined  the following relation holds:
\begin{equation}\label{eq:fundeq}
j^r\pi(j^r_D\tau_k)=j^r_D\phi_k.
\end{equation}
\end{itemize}

We will move onto explaining why the sequence of bundles ${\mathcal{J}}^r_D(M,{\mathbb C}{\mathbb P}^m)$
 are the right linearization
 of the bundles  ${\mathcal{J}}_{\rm{CR}}^r(M,{\mathbb C}{\mathbb P}^m)$, and we will also introduce the analog
of the Thom-Boardman stratification in the linearized setting.

We fix a family of so called approximately holomorphic charts
  $\varphi_{k,x}\colon ({\mathbb C}^n\times {\mathbb R},0)\rightarrow (M,x)$ (definition 3.1 in \cite{Ma04};
 see also section 3 in \cite{Ma04} for more details on their construction). This is a notion which only uses 
that $D$ is a codimension one distribution and $J$ an almost complex structure on it. We also demand
the charts to be given by CR maps as well. It is easy to see that this is always possible since our 
coordinates can be chosen to be restriction of holomorphic coordinates
 in the ambient projective space. Note that in the Levi-flat case  $D_h$ -the canonical foliation of 
${\mathbb C}^n\times \mathbb R$ by complex hyperplanes- is sent to $\mathcal{F}$.

 We treat the Levi-flat case in more detail since it is the one for which our main theorem applies:
 for each point $x\in M$, k $\gg$ 1,  and associated to the 
CR coordinates  $z^1_k,\dots,z_k^n,s_k$ defined by $\varphi_{k,x}$
 over the ball $B_{g_k}(x,\rho)$, $\rho>0$ independent of k, $x$, we have the local bundles
 ${\mathcal{J}}^r_{D_h,n,m}$ of CR r-jets with a canonical bundle map
\begin{equation}\label{eq:bmap1}
\Psi_{k,x,i}^{\rm{lin}}\colon
{\mathcal{J}}^r_D(M,{\mathbb C}^m)_i\rightarrow {\mathcal{J}}^r_{D_h,n,m}
\end{equation}
 obtained as follows:  the basis $dz_k^1,\dots, dz^n_k\in \Gamma(D^{*1,0})$ identifies $D^{*1,0}$
  with $T^{*1,0}{\mathbb C}^n$; let $I$ be an (N+2)-tuple $I=(i_0,i_1,\dots,i_n,i)$,
  $1\leq i_0\leq m$, $0\leq i_j\leq r$, $i=0,\dots,m$, $i_1+\cdots +i_n=r$. The frame
\begin{equation}\label{eqn:canframe}
  \mu_{k,x,I}:={dz_k^1}^{\odot i_1}\odot \cdots \odot {dz_k^n}^{\odot i_n}\otimes \xi_{i,i_0}
  \end{equation}
defines the bundle map of equation \ref{eq:bmap1}.

 The local bundles ${\mathcal{J}}^r_{D_h,n,m}$  glue into the non-linear
  bundle ${\mathcal{J}}^r_{\rm{CR}}(M,{\mathbb C}^m)_i$: let $y\in M$
   be a point belonging to two different charts centered at $x_0$ and $x_1$
    respectively. If we send  $y$ in both charts to the origin via a translation,
     then the change of coordinates restricts to the leaf through the origin to
     a bi-holomorphic map fixing the origin. The fibers over $y$ are  related  by
      the action of the holomorphic r-jet of the bi-holomorphism.
        If we only take the linear part of the action, there is an induced vector bundle map
\begin{equation}\label{eq:bmap2}
\Psi_{k,x_0,x_1,i}^{\rm{lin}}\colon {\mathcal{J}}^r_{D_h,n,m}\rightarrow {\mathcal{J}}^r_{D_h,n,m}
\end{equation}
which defines a vector bundle, for the cocycle condition still holds. This vector bundle
 is  ${\mathcal{J}}_{D}^r(M,{\mathbb C}^m)_i$
 as defined in   equation \ref{eqn:projjets} (it is rather a sequence of bundles
  in which the metric in the $D^{*1,0}$ factors is induced from $g_k$). Thus for  Levi-flat CR
  manifolds the vector bundles ${\mathcal{J}}_{D}^r(M,{\mathbb C}^m)_i$ are ``linear approximations'' for k $\gg$ 1
   of the non-linear bundles ${\mathcal{J}}^r_{\rm{CR}}(M,{\mathbb C}^m)_i$.  

Notice that  to make sense of equations \ref{eq:bmap1}, \ref{eqn:canframe}, and \ref{eq:bmap2}
 one does not quite need $(M,D,J)$ to be a CR manifold, just  $D$ a codimension one distribution 
endowed with an almost complex structure $J$ is enough.

The third piece of data to be linearized is the CR Thom-Boardman stratification. We have seen that 
the vector bundles  ${\mathcal{J}}_{D}^r(M,{\mathbb C}^m)_i$ and  the fiber bundles
 ${\mathcal{J}}_{\rm{CR}}^r(M,{\mathbb C}^m)_i$ use the same building blocks  ${\mathcal{J}}^r_{D_h,n,m}$,
but different transition maps. We will apply the same idea for the stratifications. Each local bundle 
${\mathcal{J}}^r_{D_h,n,m}$ carries a corresponding
  CR Thom-Boardman stratification ${\mathcal{T}}_{n,m}^{\mathcal{F}}$ (or rather a refinement which
  is a Whitney (A) stratification  \cite{Ma75}, which is invariant by r-jet extensions of bi-holomorphic
 transformations). The CR Thom-Boardman stratification
   ${\mathbb P}{\mathcal{T}}_i^{\mathcal{F}}$ of ${\mathcal{J}}^r_{\rm{CR}}(M,{\mathbb C}^m)_i$ is the
result of gluing the local stratifications
${\mathcal{T}}_{n,m}^{\mathcal{F}}$, for as mentioned these are preserved by the
transition maps; the stratifications ${\mathbb P}{\mathcal{T}}_i^{\mathcal{F}}$ in turn are used to
build   ${\mathbb P}{\mathcal{T}}^{\mathcal{F}}$ the Thom-Boardman
stratification of
${\mathcal{J}}_{\rm{CR}}^r(M,{\mathbb C}{\mathbb P}^m)$. The
relevant observation is that ${\mathcal{T}}_{n,m}^{\mathcal{F}}$ is
also preserved by $\Psi_{k,x_0,x_1,i}^{\rm{lin}}$, thus  giving
rise to the Thom-Boardman-Auroux stratifications
${\mathbb P}{\mathcal{T}}_{k,i}$ of ${\mathcal{J}}^r_D(M,{\mathbb C}^m)_i$; these stratifications in turn are
also compatible with the gluing that defines
${\mathcal{J}}^r_D(M,{\mathbb C}{\mathbb P}^m)$ giving rise to the Thom-Boardman-Auroux stratifications
 ${\mathbb P}{\mathcal{T}}_k$ (\cite{Ma04}, definition 6.4).

For general projective CR manifolds we just work with CR 1-jets. The linearization process does
 not change neither the bundle nor the 1-jet. The only thing to bear in mind is that we work with
 sequences of bundles k $\gg$ 1, for which the total space of the bundle is the same but the metric
 changes (or rather the metric defined in the vector bundles associated to the canonical affine charts of projective space).

We say that a sequence of A.H. sections $\tau_k$ of ${E_k}_{\mid M}$ is r-generic if the following conditions hold:
\begin{itemize}
 \item[(1')] $\tau_k$ is uniformly transverse along $\mathcal{F}$ to $Z_k$. 
\item[(2')] $j^r_D\phi_k$ is uniformly transverse  along $\mathcal{F}$ to ${\mathbb P}{\mathcal{T}}_k$.
\end{itemize}
We refer the reader to  \cite{Ma04}, section 4, for the notion of uniform
transversality along distributions to stratifications, but we just note that it uses a quantification of transversality, 
and requires its independence on k, for k $\gg$ 1.

\subsection{Existence of CR solutions to the linearized problem}\label{ssec:lin2}

The existence of r-generic sequences of A.H. sections of ${E_k}_{\mid M}$ is the content of 
theorem 8.4 and proposition 6.3 in \cite{Ma04}:

\begin{theorem}\label{thm:main2} Fix any
 $\delta>0$ and $r,h\in \mathbb N$, $h-r\geq 2$. Then a constant
 $\eta>0$ and a natural number $k_0$ exist such that for any
  $C^{h}$-A.H. sequence $\sigma_k$  of $E_k$ it is possible to find
   a  $C^{h}$-A.H. sequence   $s_k$ of $E_k$, so that for any $k$ bigger than   $k_0$
 \begin{itemize}
 \item
$|\nabla^j(s_k-\sigma_k)|_{g_k}<\delta,j=0,\dots,h$.
 \item Let $\tau_k$ denote the restriction of $s_k$ to $M$ and $\phi_k$ its projectivization. Then $\tau_k$ is
$\eta$-transverse along $\mathcal{F}$ to $Z_k$ and $j^r_D\phi_k$ is
$\eta$-transverse along $\mathcal{F}$ to
${\mathbb P}{\mathcal{T}}_k$.
  \end{itemize}
 \end{theorem}

Our aim in this subsection is to make sure that the solution $s_k$ in theorem \ref{thm:main2} can be 
 chosen to be  a  uniformly bounded sequence of holomorphic sections of $E_k$. For that we need to briefly recall how
theorem \ref{thm:main2} is proved: the Thom-Boardman-Auroux stratification ${\mathbb P}{\mathcal{T}}_k$ is pulled back 
to a stratification
  $j^r\pi^*{\mathbb P}{\mathcal{T}}_k\cup Z_k$ of ${\mathcal{J}}^r_D{E_k}_{\mid M}$. Theorem 7.2 in \cite{Ma04} grants the
existence of a small enough perturbation $s_k$ of $\sigma_k$, so that 
 $j^r_D\tau_k$ is uniformly  transverse to $j^r\pi^*{\mathbb P}{\mathcal{T}}_k\cup Z_k$. This encompasses both conditions
(1') and (2'): condition (1') is fulfilled by adding $Z_k$ to the stratification of ${\mathcal{J}}^r_D{E_k}_{\mid M}$. That
   $j^r_D\phi_k$ is uniformly transverse to ${\mathbb P}{\mathcal{T}}_k$ is mostly a consequence
 of equation \ref{eq:fundeq}. Then proposition
6.3 in \cite{Ma04}  implies that because the sections $\tau_k$ are A. H., uniform transversality gives as well the seemingly
 stronger uniform transversality along $\mathcal{F}$. Summarizing, the 
uniform transversality problem for $B_k$ and maps $\phi_k\colon M\backslash B_k\rightarrow {\mathbb C}{\mathbb P}^m$
 is reduced to a 
uniform transversality problem for sections $\tau_k\colon M\rightarrow {E_k}_{\mid M}$, whose solution is given as the
restriction of sections $s_k\colon {\mathbb C}{\mathbb P}^N\rightarrow E_k$. Thus we need to make sure that theorem 7.2
 in \cite{Ma04} can produce a solution $s_k$ which is a uniformly
bounded sequence of holomorphic sections of $E_k$. 

Theorem 7.2 in \cite{Ma04} gives in general a
solution of a uniform transversality problem in ${\mathcal{J}}^r_D{E_k}_{\mid M}\rightarrow M$
 by applying  relative approximately holomorphic theory as follows: in the complex manifold 
${\mathcal{N}}_\epsilon(M)\hookrightarrow  {\mathbb C}{\mathbb P}^N$ one defines as in
 equation \ref{eqn:djets}  the bundle  of pseudo-holomorphic r-jets 
\[{\mathcal{J}}^r{E_k}:=(\sum_{j=0}^r T^{*1,0}{\mathcal{N}}_\epsilon(M)\odot\cdots^{(j)}\cdots  \odot  T^{*1,0}{\mathcal{N}}_\epsilon(M))\otimes {E_k},\]
and the pseudo-holomorphic r-jet $j^r\sigma_k$ of a section $\sigma_k\colon {\mathcal{N}}_\epsilon(M)\rightarrow E_k$.
 This is done exactly as we did for $M$ but replacing the almost complex distribution $(D,J)$
 by the almost complex distribution $T{\mathcal{N}}_\epsilon(M)$.  It is also necessary to 
    ``thicken'' $j^r\pi^*{\mathbb P}{\mathcal{T}}_k\cup Z_k$ to an
appropriate stratification ${\mathcal{T}}_k$ of ${\mathcal{J}}^rE_k$. Approximately holomorphic
 theory then produces for any given A.H. sequence $\sigma_k$ of $E_k$  an arbitrary small perturbation $s_k$ so that 
$j^rs_k$ is uniformly transverse along $M$  to ${\mathcal{T}}_k$, and one checks that this implies that 
the restriction $\tau_k$ to $M$ is uniformly
 transverse to $j^r\pi^*{\mathbb P}{\mathcal{T}}_k\cup Z_k$. In other words, relative transversality theory 
associates to a suitable uniform transversality problem in the CR manifold $(M,D,J)$, a uniform transversality problem 
in the complex manifold ${\mathcal{N}}_\epsilon(M)$ so that a solution to the latter restricts to a solution to the former.

So everything is reduced to show that in a complex manifold, when one starts with a uniformly 
bounded sequence of holomorphic functions $\sigma_k$ -for example the sequence identically zero-
 the perturbations produced by approximately holomorphic theory to solve admissible uniform transversality problems
can be chosen to be holomorphic. But this is essentially the content of \cite{Do96}, section 7 (a result which is
already present in \cite{Ti89}): if $s_{k,x,i_0}^{\rm{ref}}$,
$0\leq i_0\leq m$, $x\in M$, k $\gg$ 1, is an appropriate family of
so called reference frames of $E_k=\underline{\mathbb C}^{m+1}\otimes
{\mathcal{O}}(k)$ (see section 2.3 in \cite{Au01}), then their
$L^2$-projection onto the holomorphic sections defines a family of
holomorphic reference frames. Exactly the same ideas show that for
$I=(i_1,\dots,i_N,i)$, $0\leq i\leq m$, $0\leq i_j\leq r$,
$i_1+\cdots +i_N=r$,   the $L^2$-projection of

\[\nu^{\rm{ref}}_{k,x,I}:={(z_k^1)}^{i_1}\cdots {(z_k^N)}^{i_N}s_{k,x,i}^{\rm{ref}}\]
are holomorphic sections whose pseudo-holomorphic r-jets define
reference frames of ${\mathcal{J}}^rE_k$. The perturbation provided by approximately
 holomorphic theory can be arranged to be   a finite complex linear combination of the
 holomorphic part of the $\nu^{\rm{ref}}_{k,x,I}$ (see for example the original construction in \cite{Do96}, section 3), therefore it is also holomorphic. Thus, we can 
always find a uniformly bounded sequence of holomorphic sections of $E_k$ whose restriction to $M$
is r-generic, and this proves part (B).

\subsection{Comparison between pseudo-holomorphic jets and CR jets}\label{ssec:phol}

Let $s_k$ be a uniformly bounded  sequence of holomorphic sections
 of $\underline{\mathbb C}^{m+1}\otimes {\mathcal{O}}(k)\rightarrow {\mathbb C}{\mathbb P}^N$
 provided by theorem \ref{thm:main2} (for example a perturbation
  of the trivial sequence $\sigma_k=0$). We want to  check part (C) of our strategy.  

The restriction of $s_k$ to $M$ is a sequence $\tau_k$ of CR sections. Because $\tau_k$ is 
transverse to $Z_k$ along $\mathcal{F}$ for k $\gg$ 1, its zero set $B_k$ is a CR 
 submanifold of the expected dimension, and this proves point (1) 
in theorem \ref{thm:gener}.

 Let $\phi_k\colon M\backslash B_k\rightarrow {\mathbb C}{\mathbb P}^m$ be the corresponding
 sequence of CR maps.  By hypothesis $j^r_D\phi_k\in \Gamma({\mathcal{J}}^r_D(M\backslash B_k,{\mathbb C}{\mathbb P}^m))$
  is uniformly transverse along $\mathcal{F}$ to ${\mathbb P}{\mathcal{T}}_k$ for k $\gg$ 1. 
The CR r-jet  $j^r_{\rm{CR}}\phi_k$ is a section of
${\mathcal{J}}^r_{\rm{CR}}(M\backslash
B_k,{\mathbb C}{\mathbb P}^m)$ (the same bundle for all k, apart
from the submanifold of base points). The key observation is again that  for k $\gg$ 1 the pairs
$({\mathcal{J}}^r_D(M,{\mathbb C}{\mathbb P}^m),{\mathbb P}{\mathcal{T}}_k)$
 and $({\mathcal{J}}^r_{\rm{CR}}(M,{\mathbb C}{\mathbb P}^m),{\mathbb P}{\mathcal{T}}^{\mathcal{F}})$ are constructed
 using the same building blocks $({\mathcal{J}}^r_{D_h,n,m},{\mathcal{T}}_{n,m}^{\mathcal{F}})$,
 but with different transition functions.

 For each $x\in M\backslash B_k$, fix $i$ such  that
 $j^r_D\phi_k\colon B_{g_k}(x,\rho)\backslash B_k\rightarrow  {\mathcal{J}}^r_D(M\backslash B_k,{\mathbb C}^m)_i$.
  Over $B_{g_k}(x,\rho)$ we use the bundle map of equation \ref{eq:bmap1}
   to see $j^r_D\phi_k$ as a section of ${\mathcal{J}}^r_{D_h,n,m}$. Similarly, we can
    see $j^r_{\rm{CR}}\phi_k$ over $B_{g_k}(x,\rho)\backslash B_k$ as a section of ${\mathcal{J}}^r_{D_h,n,m}$.

The checking we have to make requires different ideas depending on whether we stay uniformly away from $B_k$ -i.e. at some small
 fixed distance independently of k- or not.

Because $\phi_k$ is a projectivization if we stay uniformly away from $B_k$  then $\phi_k$ is a sequence of 
uniformly bounded CR sections
 (for either $g_k$ or the local standard metric $g_0$ in $B_{g_k}(x,\rho)\backslash B_k$, the flat connection and the standard
 Hermitian metric). In such a situation on easily checks

\begin{equation}\label{eq:distance}
|j^r_D\phi_k-j^r_{\rm{CR}}\phi_k|_{C^1,g_k}\leq O(k^{-1/2}),
\end{equation}
for all k $\gg$ 1.

Uniform transversality along $\mathcal{F}$ is an open condition (see
\cite{Ma04}, section 7),
 meaning that if $j^r_D\phi_k$ is  $\delta$-transverse along $\mathcal{F}$ to ${\mathbb P}{\mathcal{T}}_k$
  and $|j^r_D\phi_k-\xi_k|_{C^1}\ll \delta$, then $\xi_k$ is $\delta/2$-transverse
  along $\mathcal{F}$ to ${\mathbb P}{\mathcal{T}}_k$. But this is exactly what equation \ref{eq:distance}
 grants in the points uniformly away from $B_k$. The conclusion is that on each 
$B_{g_k}(x,\rho)\backslash U_{k,x,\epsilon}$, where  $U_{k,x,\epsilon}$ is
 the neighborhood of $B_k\cap B_{g_k}(x,\rho)$ of uniform radius $\epsilon$ independent of $x$, k, the local section  
$j^r_{\rm{CR}}\phi_k\in \Gamma({\mathcal{J}}^r_{D_h,n,m})$ is uniformly transverse to the linearized
 Thom-Boardman stratification ${\mathcal{T}}_{n,m}^{\mathcal{F}}$, which is also the CR Thom-Boardman
 stratification.

As for the points in $U_{k,x,\epsilon}$,  due to the definition of $j^r\pi^*{\mathbb P}{\mathcal{T}}_k\cup Z_k$
 the graph of $j^r_D\tau_k$  stays uniformly away from the points in $Z_k$ to which the other strata 
 of $j^r\pi^*{\mathbb P}{\mathcal{T}}_k\cup Z_k$ converge, and this is something that just have
 to do with the  1-jet prolongation $j^1_D\tau_k$. This implies that in $U_{k,x,\epsilon}$  $j^1_D\phi_k$ stays uniformly away from all the strata of
 ${\mathbb P}{\mathcal{T}}_k$ of strictly positive codimension  (\cite{Ma04}, remark 6.2). By the construction of 
the Thom-Boardman-Auroux stratification any section of ${\mathcal{J}}^r_D(M,{\mathbb C}{\mathbb P}^m)$
 whose degree 1 part coincides with  $j^1_D\phi_k$ will also stay uniformly away from these 
strata. Thus inside $U_{k,x,\epsilon}$   $j^r_{\rm{CR}}\phi_k\in \Gamma({\mathcal{J}}^r_{D_h,n,m})$ -whose degree 1 part is
 $j^1_D\phi_k\in \Gamma({\mathcal{J}}^1_{D_h,n,m})$- will also stay uniformly away from
 the strata of strictly positive codimension. Therefore by its very definition uniform transversality
 for $j^r_{\rm{CR}}\phi_k$ to ${\mathbb P}{\mathcal{T}}^{\mathcal{F}}$ holds as well in $U_{k,x,\epsilon}$, and this proves
part (C).

\subsection{Proof of theorem \ref{thm:gener}} Given $h\in \mathbb{N}$ we apply theorem \ref{thm:OS} to embed $(M,D,J)$ as a 
$C^h$-Levi-flat submanifold of some projective space. Then the previous three subsections imply that we
 can always construct $s_k$ uniformly bounded sequences of holomorphic sections of 
$\underline{\mathbb C}^{m+1}\otimes {\mathcal{O}}(k)$,
  such that (1) the set of base points $B_k$ of the restriction to $M$  is a CR submanifold of codimension 2m+2 and class
$C^h$, and (2)  the projectivization of the restriction  $\phi_k\colon M\backslash B_k\rightarrow {\mathbb C}{\mathbb P}^m$
is a  r-generic CR map.

To prove point (3) in the statement of the theorem \ref{thm:gener} we have to
 solve yet another uniform transversality problem:
one considers the bundles
${\mathcal{J}}^1_{TM^J}({\mathcal{N}}_\epsilon(M),{\mathbb C}{\mathbb P}^m)$
of pseudo-holomorphic 1-jets along $TM^J$  for maps to
${\mathbb C}{\mathbb P}^m$, defined in section \ref{sec:sec1}.
We can equally define them using affine coordinates of
${\mathbb C}{\mathbb P}^m$ and then gluing the  vector
bundles
${\mathcal{J}}^1_{TM^J}({\mathcal{N}}_\epsilon(M),{\mathbb C}^m)_i$.
These are sequences of vector bundles with metric depending on k.
Inside we have the submanifolds $\Sigma_{r,k}$ of 1-jets along $TM^J$
whose degree 1 homogeneous component has rank $r\in \{0,\dots,m\}$.  It is possible to pull them back
to strata $\Sigma_{r,k}^J$ of the bundles ${\mathcal{J}}^1E_k\rightarrow
{\mathcal{N}}_\epsilon(M)$ of pseudo-holomorphic 1-jets of $E_k$.

The stratification $Z_k\cup \Sigma_{r,k}^{J}$ of ${\mathcal{J}}^1E_k$ is
such that there is a version of  theorem \ref{thm:main2} asserting that for any $\delta'>0$ there exist
$\eta'>0$ and  $s_k'$ a uniformly bounded sequence of holomorphic
section of $E_k$, such that for all k $\gg$ 1
\begin{itemize}
\item $|\nabla^j(s_k-s_k')|_{g_k}<\delta',j=0,\dots,h$
\item $j^1_{TM^j}\Phi_k'$ is $\eta'$-transverse along $M$ to $\Sigma_{r,k}$, where
 $\Phi_k'$ is the projectivization of $s_k'$.
\end{itemize}

By choosing $\delta''$ small enough the sequence $s_k'$ will be such that $\phi'_k$ will still have properties
 (1) and (2) at the beginning of this subsection, and this proves theorem \ref{thm:gener}.

 \subsection{Existence of Lefschetz pencil structures and proofs of 
theorems \ref{thm:pencil} and \ref{thm:hyp}}

In the Lefschetz pencil condition we ask the holomorphic sections $s_k$ to have components with zero locus intersecting
in $\tilde{B_k}$ a real codimension 4 holomorphic submanifold, which intersects $M$ transversaly in $B_k$. 
This is yet another transversality problem for 0-jets. Namely, one asks that each summand of $s_k$  intersects the zero section of 
the corresponding summand of $E_k$ transversaly, that $s_k$ intersects the zero section of $E_k$ transversaly, and that  
that $s_k$ intersects the zero section of $E_k$ transversaly along $M$.

Property (2) it attained with a small perturbation following mostly from the proof  of
theorem \ref{thm:gener}. Indeed, as we mentioned at the end of subsection \ref{ssec:lin1}, for 1-jets
 the linearized version of the CR jets coincides
with the CR 1-jet as introduced in definition \ref{def:1jets}. 
The linearized version of the bundles of CR 1-jets is the same bundle of CR
 1-jets, but having in mind that in the sequence the only piece of data that varies is the metric.

Subsection \ref{ssec:lin2} is an account in first place  how to solve the r-genericity problem in the framework of approximately
holomorphic geometry, and no assumption on the integrability of either $D$ or $J$ is needed. Here we just use the result
granting the existence of sequences of 1-generic A.H. sections of $E_k$.  Then there is a sketch of how to get solutions
which are uniformly bounded sequences of holomorphic sections. The only necessary ingredient is the fact that the ambient symplectic
manifold admits a compatible complex structure. 

Subsection \ref{ssec:phol} explains how for a Levi-flat CR manifold, the pseudo-holomorphic r-jets and the CR r-jets are as
close as desired as $k$ grows. Since for 1-jets there is no difference between both concepts, we do not need to use it here and
we can conclude the existence holomorphic sections whose restriction to our arbitrary projective CR manifold is 1-generic. 

Point (3) in theorem \ref{thm:gener} coincides with condition (3) in the definition of Lefschetz pencil, so we can proceed
as we did in the proof of theorem \ref{thm:gener} to get it (which again requires no integrability of $D$).

The part concerning the ``Lefschetz hyperplane theorem'' for homotopy
 (homology) groups of a hyperplane section is a standard result valid
  for Lefschetz pencil structures for 2-calibrated structures (\cite{IM04a}, corollary 1.2), and this finishes the proof 
of theorem \ref{thm:pencil}.

Back to Levi-flat CR manidfolds, proposition 1 is also elementary: coordinates about $B_k$ are again constructed
 as for  Lefschetz pencil structures for 2-calibrated structures (\cite{IM04}, theorem 1.2), but we obtain CR 
coordinates because the two components of
  $\tau_k\in \Gamma(\underline{{\mathbb C}^2}\otimes {\mathcal{O}(k)}_{\mid M})$ are already CR;
  CR Morse coordinates about $\Delta_k$ are obtained by applying the complex Morse lemma with parameters.

Theorem \ref{thm:hyp} is proved in \cite{Ma05} for general 2-calibrated foliations and for Lefschetz 
pencils so that the normal forms of proposition \ref{pro:nforms} hold and the curve $\phi_k(\Delta)$
  is an immersed curve with generic self-intersections. One observes that the proof easily generalizes
 for $\phi_k(\Delta)$  an immersed curve, which is what theorem \ref{thm:pencil} grants.


\begin{thebibliography}{xxxxx}

\bibitem{Au97} D. Auroux, Asymptotically holomorphic families of symplectic submanifolds. Geom.
Funct. Anal. {\bf 7} (1997), 971-995.


\bibitem{Au00}  D. Auroux , Symplectic 4-manifolds as branched coverings of ${\mathbb C}{\mathbb P}^2$.
 Invent. Math. {\bf 139} (2000), 551-602.

\bibitem{Au01} D. Auroux, Estimated transversality in symplectic geometry and projective maps.
  Symplectic geometry and mirror symmetry (Seoul, 2000),  1-30, World Sci. Publishing, River Edge, NJ, 2001.




\bibitem{AK00} D. Auroux, L. Katzarkov,  Branched coverings of ${\mathbb C}{\mathbb P}^2$  
and invariants of symplectic 4-manifolds.
Inventiones Math. {\bf 142} (2000), 631-673.

\bibitem{Be02} M. Bertelson, A h-principle for open relations invariant under foliated isotopies.
 J. Symplectic Geom.  {\bf 1} (2002), no. 2, 369-425.

\bibitem{Bo67} J. M. Boardman , Singularities of differentiable maps.   Publ. Math.  Inst. Hautes
 \'Etudes Sci.  {\bf 33} (1967), 21-57.

\bibitem{BdM75} L. Boutet de Monvel,  Int\'egration des \'equations de Cauchy-Riemann induites formelles.
  S\'eminaire Goulaouic-Lions-Schwartz 1974-1975; \'Equations aux deriv\'ees partielles lin\'eaires 
et non lin\'eaires,  pp. Exp. {\bf 9}, 14 pp. Centre Math., \'Ecole Polytech., Paris (1975).


\bibitem{CSW04} J. Cao, M. Shaw, L. Wang, Estimates for the d-bar-Neumann problem and nonexistence
 of $C^2$ Levi-flat hypersurfaces in  ${\mathbb C}{\mathbb P}^n$.  Math. Zeit.  {\bf 248} (2004), 183-221.

\bibitem{CSW05} J. Cao, M. Shaw, The $\overline\partial$-Cauchy problem and nonexistence 
of Lipschitz Levi-flat hypersurfaces in ${\mathbb C}{\mathbb P}^n$ with $n\geq 3$.  Math. Z. {\bf 256}  (2007),  no. 1, 175–192,


\bibitem{SGA7II} P. Deligne, N.  Katz,  S\`eminaire de G\`eom\`etrie Alg\'ebrique du Bois-Marie, 1967-1969.
Lecture Notes in Math. {\bf 340} (1973).

\bibitem{De08} B. Deroin,  Laminations dans les espaces projectifs complexes. (French) [Laminations
 in complex projective spaces]  J. Inst. Math. Jussieu {\bf 7}  (2008),  no. 1, 67-91

\bibitem{Do96} S. K. Donaldson, Symplectic submanifolds and almost-complex geometry. J. Differential
Geom. {\bf 44} (1996), 666-705.

\bibitem{Do98} S. K. Donaldson, Lefschetz fibrations in symplectic geometry. Doc. Math. Extra Vol.
 ICM 98, {\bf II} (1998), 309-314.

\bibitem{DS00} S.K. Donaldson, I. Smith, Lefschetz pencils and the canonical class for symplectic 
four-manifolds. Topology {\bf 42} (2003),  no. 4, 743-785.

\bibitem{Gh99} E. Ghys,  Laminations par surfaces de Riemann. Dynamique et g\'eom\'etrie complexes
 (Lyon, 1997), ix, xi, 49-95, Panor. Synth\`eses {\bf 8}, Soc. Math. France, Paris (1999).

\bibitem{Gi02} E. Giroux, G\'eom\'etrie de contact de la dimension trois vers les dimensions sup\'erior,
Proceedings of the ICM, Beijing 2002, vol. {\bf 2}, 405-414.

\bibitem{GM01} E. Giroux and J.-P. Mohsen, Structures de contact et fibrations symplectiques
sur le cercle.  In preparation.

\bibitem{Ho}  L. H\"ormander,   Notions of convexity. Progress in Mathematics {\bf 127}, 
Birkh\"auser Boston, Inc., Boston, MA (1994).

\bibitem{IM04a} A. Ibort, D. Mart\'{\i}nez Torres,  Approximately holomorphic geometry and
 estimated transversality on 2-calibrated manifolds.  C. R. Math. Acad. Sci. Paris  {\bf 338}  (2004),  no. 9, 709–712.

\bibitem{IM04}  A. Ibort, D. Mart\'{\i}nez Torres, Lefschetz pencil structures for 2-calibrated manifolds.
 Comptes Rendus Mathematique  {\bf 339} (2004), Issue 3, 215-218.



\bibitem{MY04} G. Marinescu, N. Yeganefar, Embeddability of some strongly pseudoconvex CR manifolds.
  Trans. Amer. Math. Soc.  {\bf 359}  (2007),  no. 10, 4757--4771


\bibitem{Ma04} D. Mart\'{\i}nez Torres,  The geometry of 2-calibrated manifolds.  Port. Math.  {\bf 66}  (2009),  no. 4, 427-512

\bibitem{Ma05}   D. Martínez Torres, Codimension-one foliations calibrated by nondegenerate closed 2-forms. Pacific J. Math. {\bf 261} (2013), no. 1, 165--217

\bibitem{Ma08} D. Martínez Torres, A note on strict $\mathbb{C}$-convexity. Rev. Mat. Complut. {\bf 25} (2012), no. 1, 125--137.

\bibitem{Ma75} J. N. Mather, How to stratify mappings and jet spaces. Singularit\'es 
d'Applications Diff\'erentiables (Plans-sur-Bex 1975), Lecture Notes in Math {\bf 535}, Springer (1976), 128-176.

\bibitem{MPS02}  V. Mu\~noz, F. Presas, I. Sols, Almost holomorphic embeddings in Grassmannians
 with applications to singular symplectic submanifolds. J. Reine Angew. Math.  {\bf 547}  (2002), 149-189.


\bibitem{NW03} L. Ni, J. Wolfson,  The Lefschetz theorem for CR submanifolds and the nonexistence of real analytic
Levi flat submanifolds. Comm. in Anal. and Geom. {\bf 11} (2003), 553-564.


\bibitem{OS00} T. Ohsawa, N. Sibony,  K\"ahler identity on Levi flat manifolds and application to
 the embedding. Nagoya Math. J. {\bf 158} (2000), 87-93.


\bibitem{Su76} D. Sullivan, Cycles for the dynamical study of foliated manifolds and complex manifolds.
  Invent. Math. {\bf 36}  (1976), 225-255.

\bibitem{Ti89}  G. Tian, On a set of polarized K\"ahler metrics on algebraic manifolds. 
J. Differential Geom.  {\bf 32} (1990), no. 1, 99-130.














\end{thebibliography}
\end{document}